\documentclass{article}
\usepackage{maa-monthly}

\theoremstyle{theorem}
\newtheorem{theorem}{Theorem}

\newtheorem{corollary}[theorem]{Corollary}

\theoremstyle{definition}
\newtheorem{definition}[theorem]{Definition}
\newtheorem{remark}[theorem]{Remark}

\begin{document}

\title{Extension of the Fundamental Theorem of Algebra for Polynomial Matrix Equations with Circulant Matrices }
\markright{Extension of the Fundamental Theorem of Algebra}
\author{Vyacheslav M. Abramov}  

\maketitle

\begin{abstract}
We establish an analogue of the fundamental theorem of algebra for polynomial matrix equations, in which the matrices-coefficients and unknown matrix are assumed to be circulant matrices.
\end{abstract}

%

\section{Introduction.}\label{S1}
The fundamental theorem of algebra (FTA) has a long and distinguished history going back to the seventeenth century, where the problem was mentioned
by Peter Roth in his book \textit{Arithmetica Philosophica} published in 1608 and by Albert Girard in his book \textit{L'invention nouvelle en l'Alg\`{e}bre} published in 1629
(see \cite{Gardner, Man}). In the eighteenth century, many attempts to prove it were due to Jean-Baptiste le Rond d'Alembert, Leonhard Euler, Fran\c{c}ois Daviet de Foncenex, Carl Friedrich Gauss, Joseph-Louis Lagrange, Pierre-Simon de Laplace, and James Wood. Carl Friedrich Gauss is often considered as a mathematician given credit for providing the first correct proof of the FTA in his 1799 doctoral dissertation. His proof, however, contained a gap that was fixed in an elementary way in \cite{BV}. Another incomplete proof based on an original idea was due to James Wood in 1798, one year before Gauss's proof. The complete proof of the theorem based on Wood's idea is provided in \cite{S}. The first textbook containing a full proof of the FTA is a book by Augustin-Louis Cauchy \cite{C}.
Nowadays there are a large variety of different proofs (e.g., \cite{FR, TU}), the simplest of which seem to be given in \cite{A, B}. Another elementary proof based on the only four basic arithmetical operations  has been recently provided in \cite{O}. The FTA for an algebraically closed field with characteristic zero is proved in \cite{Sh} (see also \cite{Al}).

In the present note, we establish an analogue of the FTA for polynomial matrix equations, in which the matrices-coefficients and unknown matrix are assumed to be circulant matrices (CM) with complex entries or, more generally, with entries belonging to an algebraically closed field with characteristic zero. Recall that a CM has the form
\[
A=\left(\begin{matrix}a_0 &a_1 &\cdots &a_{d-1}\\
a_{d-1} &a_0 &\cdots &a_{d-2}\\
\vdots &\vdots &\cdots &\vdots\\
a_1 &a_2 &\cdots &a_0
\end{matrix}\right),
\]
denoted further by $\texttt{circ}(a_0, a_1,\ldots, a_{d-1})$.

We consider the polynomial matrix equation
\begin{equation}\label{1}
\boldsymbol{X}^n+\boldsymbol{A}_1\boldsymbol{X}^{n-1}+\boldsymbol{A}_2\boldsymbol{X}^{n-2}+\ldots+\boldsymbol{A}_{n-1}\boldsymbol{X}+\boldsymbol{A}_n=\boldsymbol{O},
\end{equation}
in which all the matrices that appear on the left-hand side of the equation are assumed to be $d\times d$ CM; $\boldsymbol{O}$ denotes the $d\times d$ matrix of zeros.

The study of the classes of equation \eqref{1}, where the matrices $\boldsymbol{A}_1$, $\boldsymbol{A}_2$,\ldots, $\boldsymbol{A}_n$ and solutions $\boldsymbol{X}$ belong to the class of CM is natural. CM form a commutative ring \cite[p. 35]{G} and are attractive for the study of the matrix equations. CM have many interesting properties and play a significant role in a number of applications. For a recent study of CM and their new interesting properties see \cite{KS}. In \cite{KS, KW}, CM were used for the solution of usual polynomial equations of degrees $3$ and $4$.

\medskip
\noindent
\begin{definition}
We say that a class of polynomial matrix equations satisfies the FTA, if any polynomial matrix equation of degree $n\geq 1$ of that class has at least one solution.
\end{definition}

Let $\boldsymbol{A}_k=\texttt{circ}(a_{k,0}, a_{k,1},\ldots, a_{k,d-1})$, let $r_0$, $r_1$, $r_{2}$,\ldots, $r_{d-1}$ denote primitives of the $d$th root of unity, $r_k=\mathrm{e}^{\boldsymbol{i}2\pi k/d}$, $\boldsymbol{i}=\sqrt{-1}$, and let $\overline{r}_k=1/r_k$, $k=0,1,\ldots,d-1$ denote the conjugate of $r_k$. Assuming also that the index parameter $k$ for $r_k$ (or $\overline{r}_k$) can be greater than $d-1$ or less than $0$, we follow the convention that $r_k=r_{k(\mathrm{mod}\ d)}$ ($\overline{r}_k=\overline{r}_{k(\mathrm{mod}\ d)}$).

The main results of this note are given by the following theorems.

\begin{theorem}
Equation \eqref{1} satisfies the FTA with the total number of solutions not exceeding $n^d$.
\end{theorem}

\begin{theorem}
Let $n_i$ denote the number of distinct roots of the monic polynomial equation
\begin{align}\label{3}
x^n+b_{1}^{(i)}x^{n-1}+\ldots+b_{n}^{(i)}=0, \quad i=1,2,\ldots,d,
\end{align}
where
\[
b_{k}^{(i)}=\sum_{j=1}^{d}a_{k,j-1}\overline{r}_{(i-1)(j-1)}.
\]
Then the total number of solutions of \eqref{1} is $\prod_{i=1}^{d}n_i$.
\end{theorem}

\begin{corollary}
The total number of solutions of \eqref{1} attains $n^d$ if and only if all the roots of each of the monic polynomial equations \eqref{3} are distinct.
\end{corollary}

\medskip
The rest of the paper is organized as follows. In Section \ref{S2}, we prove Theorem 2. In Section \ref{S3}, we prove Theorem 3.

\section{Proof of Theorem 2.}\label{S2}
\subsection{Background.}
Recall the following definition. A matrix is circulant if and only if it is a linear combination of the powers of the matrix
\[
\boldsymbol{C}=\boldsymbol{C}_d=\left(\begin{matrix}0 &1 &0 &\cdots &0\\
0 &0 &1 &\cdots &0\\
\vdots &\vdots &\vdots &\ddots &\vdots\\
0 &0 &0 &\cdots &1\\
1 &0 &0 &\cdots &0\end{matrix}\right).
\]
The matrix $\boldsymbol{C}_d$ is known to be a generator of a cyclic group of order $d$ with the property $\boldsymbol{C}_d^d=\boldsymbol{I}$, where $\boldsymbol{I}$ is the identity matrix \cite{W}.   For example, for the powers of the matrix in dimension $3\times3$ we have:
\[
\boldsymbol{C}_3=\left(\begin{matrix}
0 &1 &0\\
0 &0 &1\\
1 &0 &0
\end{matrix}\right), \ \boldsymbol{C}_3^2=\left(\begin{matrix}
0 &0 &1\\
1 &0 &0\\
0 &1 &0
\end{matrix}\right), \ \boldsymbol{C}_3^3=\left(\begin{matrix}
1 &0 &0\\
0 &1 &0\\
0 &0 &1
\end{matrix}\right), \ \boldsymbol{C}_3^4=\boldsymbol{C}_3, \ \text{etc.}
\]

The matrix $\boldsymbol{C}$ is a $d\times d$ permutation matrix, hence orthogonal/unitary and hence normal \cite[p. 26]{Ar}. Thus by the spectral
theorem, it has an orthonormal basis of eigenvectors. The eigenvectors of the matrix are precisely the columns of the following matrix
\[
\boldsymbol{S}=\left(\begin{matrix}1 &1 &1 &1 &\ldots &1 \\
1 &r_1 &r_2 &r_3 &\ldots &r_{d-1} \\
1 &r_{2} &r_{4} &r_{6} &\ldots &r_{2d-2} \\
1 &r_{3} &r_{6} &r_9 &\ldots &r_{3d-3}\\
\vdots &\vdots &\vdots &\vdots &\ldots &\vdots\\
1 &r_{d-1} &r_{2d-2} &r_{3d-3} &\ldots &r_{(d-1)^2}
\end{matrix}\right),
\]
(see \cite{KS})
which satisfies the following properties: it is a symmetric Vandermonde matrix (i.e., a double Vandermonde matrix in rows and columns), and being multiplied by the factor $(\sqrt{d}/d)$ becomes a discrete Fourier transform (DFT) matrix \cite[Chapter 2]{RY} and unitary matrix.

The following proof of the theorem is based on the background of the DFT matrix.

\subsection{Proof.} The key point is that the entire algebra $CM_d$ of circulant matrices is the span
of powers of $\boldsymbol{C}$, with the permutation $\boldsymbol{C}$ having order $d$. Moreover, the first $d$ powers of $\boldsymbol{C}$ are
linearly independent. Thus, we have a natural surjection of $d$-dimensional $\mathbb{C}$-algebras (if we assume that the entries of the CM belong to the field $\mathbb{C}$)
\[
\mathbb{C}[t]/(t^d-1)\twoheadrightarrow CM_d, \quad t\mapsto \boldsymbol{C},
\]
and this must be therefore an isomorphism.

Moreover, $\boldsymbol{S}\boldsymbol{C}\boldsymbol{S}^{-1}$
is a diagonal matrix from above. Therefore so is $\boldsymbol{S}p(\boldsymbol{C})\boldsymbol{S}^{-1}=p(\boldsymbol{S}\boldsymbol{C}\boldsymbol{S}^{-1})$
for any polynomial $p(t)$. In other words, the DFT matrix simultaneously diagonalizes the entire
algebra $CM_d$.

Now suppose we have the circulant polynomial equation \eqref{1} with all the terms being $d\times d$ complex circulant matrices. 
By the first part of the proof, conjugating by $\boldsymbol{S}$ yields a polynomial equation
\begin{align}\label{5}
\boldsymbol{U}^n+\boldsymbol{B}_1\boldsymbol{U}^{n-1}+\boldsymbol{B}_2\boldsymbol{U}^{n-2}+\ldots+\boldsymbol{B}_{n-1}\boldsymbol{U}+\boldsymbol{B}_n=\boldsymbol{O},
\end{align}
where all matrices are now $d\times d$ complex diagonal matrices.

Let the diagonal entries of $\boldsymbol{U}$ be given by $u_1$, $u_2$,\ldots, $u_d$. Then the above matrix equation reduces,
diagonal-entry by diagonal-entry, to $p_i(u_i)=0$ for some $n$th degree monic polynomials $p_i\in\mathbb{C}[t]$. By the FTA each of the equations $p_i(u_i)=0$ has a solution, and by elementary field theory, $u_i$ can take at most $n$ distinct values, so that the diagonal matrix
$\boldsymbol{U}=\boldsymbol{S}\boldsymbol{X}\boldsymbol{S}^{-1}$
has at least one and at most $n^d$ possibilities. But then so does $\boldsymbol{X}$.

\begin{remark}
In fact, if we fix a size $d$, then the above proof works over any algebraically closed field with characteristic zero containing a primitive $d$th
root of unity. This is because then one avoids using the spectral theorem and directly works
with the invertible Vandermonde matrix $\boldsymbol{S}$, whose columns form an eigenbasis for $\boldsymbol{C}$, and hence a
simultaneous eigenbasis for all of $CM_d$.
\end{remark}

\section{Proof of theorem 3.}\label{S3}
Denote the entries of the matrix $\boldsymbol{S}$ by $[s_{i,j}]$, and the entries of the matrix $\boldsymbol{S}^{-1}$ by $[\tilde{s}_{i,j}]$. Then we have $s_{i,j}=r_{(i-1)(j-1)}$ and $\tilde{s}_{i,j}=d^{-1}\overline{r}_{(i-1)(j-1)}$, \ $i,j=1,2,\ldots,d$.

Let $\boldsymbol{A}=\texttt{circ}(a_0, a_1,\ldots,a_{d-1})$, and for any integer $k$, $a_k=a_{k(\textrm{mod}~d)}$. Taking into account that the entries of the matrix $\boldsymbol{A}$ are $[a_{j-i}]$,
for the diagonal entries of the matrix $\boldsymbol{V}=\boldsymbol{S}\boldsymbol{A}\boldsymbol{S}^{-1}$ denoted by $v_1$, $v_2$,\ldots, $v_d$, we have the following presentation
\[
v_i=\sum_{j=1}^{d}\sum_{k=1}^{d}s_{i,k}a_{j-k}\tilde{s}_{j,i}.
\]
Substituting $s_{i,k}$ and $\tilde{s}_{j,i}$ for their corresponding values $r_{(i-1)(k-1)}$ and $d^{-1}\overline{r}_{(j-1)(i-1)}$, we obtain:

\begin{equation}\label{7}
\begin{aligned}
v_i&=\frac{1}{d}\sum_{j=1}^{d}\sum_{k=1}^{d}a_{j-k}r_{(i-1)(k-1)}\overline{r}_{(j-1)(i-1)}\\
&=\frac{1}{d}\sum_{j=1}^{d}\sum_{k=1}^{d}a_{j-k}\overline{r}_{(i-1)(1-k)}\overline{r}_{(j-1)(i-1)}\quad \text{\bigg\{since} \ r_{(i-1)(k-1)}=\overline{r}_{(i-1)(1-k)}\bigg\}\\
&=\frac{1}{d}\sum_{j=1}^{d}\sum_{k=1}^{d}a_{j-k}\overline{r}_{(i-1)(j-k)}\quad\quad\quad\quad\quad\ \ \text{\bigg\{since} \ \overline{r}_{l}\overline{r}_{m}=\overline{r}_{l+m}\bigg\}\\
&=\sum_{j=1}^{d}a_{j-1}\overline{r}_{(i-1)(j-1)}.\quad\quad\quad\quad\quad\quad\quad\quad \text{\bigg\{by simple algebra\bigg\}}
\end{aligned}
\end{equation}

Let us now turn to matrix equation \eqref{5}. Since the columns of the matrix $\boldsymbol{S}^{-1}$ are linearly independent, then, due to presentation \eqref{7}, the change of the original variables to the new ones when we pass from equation \eqref{1} to equation \eqref{5} is indeed lawful. That is, matrix equation \eqref{5} reduces to the system of independent monic polynomial equations.

Denote the $i$th diagonal entry of the matrix $\boldsymbol{B}_k$  by $b_{k}^{(i)}$. From presentation \eqref{7} we have
\begin{align}\label{9}
b_{k}^{(i)}=\sum_{j=1}^{d}a_{k,j-1}\overline{r}_{(i-1)(j-1)}, \quad k=1,2,\ldots,n, \quad i=1,2,\ldots, d.
\end{align}
Then, the $i$th monic polynomial equation $p_i(u_i)=0$ has the form
\begin{align}\label{11}
u_i^n+b_{1}^{(i)}u_i^{n-1}+\ldots+b_{n}^{(i)}=0,
\end{align}
where $b_{k}^{(i)}$ is given by \eqref{9}. Assume now that the $i$th equation in \eqref{11} has $n_i$ distinct roots. Then the total number of solutions of \eqref{5} must be $\prod_{i=1}^{d}n_i$. It is also true for the total number of solutions of \eqref{1}. Consequently, if each of the equations has distinct roots, i.e., $n_i\equiv n$, then the total number of solutions of \eqref{1} is $n^d$.

\begin{acknowledgment}{Acknowledgments.}
The author sincerely thanks the referees for their comprehensive reviews.
\end{acknowledgment}

\begin{affil}
24 Sagan Drive, Cranbourne North, Melbourne, Victoria-3977, Australia\\
vabramov126@gmail.com
\end{affil}

\vfill\eject

\end{document}